\documentclass[leqno]{article}
\usepackage{amsmath}
\usepackage{amssymb}
\newcommand{\F}{{\mathbb{F}}}
\newcommand{\C}{{\mathbb{C}}}
\newcommand{\Z}{{\mathbb{Z}}}
\newcommand{\ba}{{\boldsymbol{a}}}
\newcommand{\bb}{{\boldsymbol{b}}}
\newcommand{\bc}{{\boldsymbol{c}}}
\newcommand{\bH}{{\mathbf{H}}}
\newcommand{\fC}{{\mathfrak{C}}}
\newcommand{\fI}{{\mathfrak{I}}}
\newcommand{\fS}{{\mathfrak{S}}}
\newcommand{\fX}{{\mathfrak{X}}}
\newcommand{\cB}{{\mathcal{B}}}
\newcommand{\cC}{{\mathcal{C}}}
\newcommand{\cF}{{\mathcal{F}}}
\newcommand{\cL}{{\mathcal{L}}}
\newcommand{\cR}{{\mathcal{R}}}
\newcommand{\cLR}{{\mathcal{LR}}}
\newcommand{\cO}{{\mathcal{O}}}
\newcommand{\cS}{{\mathcal{S}}}
\newcommand{\cU}{{\mathcal{U}}}
\newcommand\Irr{\operatorname{Irr}}
\newcommand\Ind{\operatorname{Ind}}
\newcommand\sgn{\varepsilon}
\renewcommand{\leq}{\leqslant}
\renewcommand{\geq}{\geqslant}

\newtheorem{thm}{Theorem}[section]
\newtheorem{cor}[thm]{Corollary}
\newtheorem{prop}[thm]{Proposition}
\newtheorem{lem}[thm]{Lemma}
\newtheorem{conj}[thm]{Conjecture}

\newenvironment{rem}{\refstepcounter{thm}
\medskip \noindent {\em Remark \thethm} }{\par\medskip}

\newenvironment{defn}{\refstepcounter{thm}
\medskip \noindent {\bf Definition \thethm} }{\par\medskip}

\newenvironment{exmp}{\refstepcounter{thm}
\medskip \noindent {\bf Example \thethm} }{\par\medskip}

\newenvironment{proof}{\noindent {\it Proof.}}{\hfill $\Box$ \smallskip}

\begin{document}

\date{}

\title{On the Kazhdan--Lusztig order on cells and families}

\author{Meinolf Geck \\ {\footnotesize Institute of Mathematics, King's 
College, Aberdeen AB24 3UE, Scotland, UK} \\ {\footnotesize\it E-mail 
address: \tt m.geck@abdn.ac.uk}}

\maketitle
\pagestyle{myheadings}
\markboth{Geck}{Kazhdan--Lusztig order on cells and families}

\begin{abstract} 
We consider the set $\Irr(W)$ of (complex) irreducible characters of a 
finite Coxeter group~$W$. The Kazhdan--Lusztig theory of cells gives rise
to a partition of $\Irr(W)$ into ``families'' and to a natural partial
order $\leq_{\cLR}$ on these families. Following an idea of Spaltenstein, 
we show that $\leq_{\cLR}$ can be characterised (and effectively computed) 
in terms of standard operations in the character ring of $W$. If, moreover,
$W$ is the Weyl group of an algebraic group $G$, then $\leq_{\cLR}$ can 
be interpreted, via the Springer correspondence, in terms of the closure 
relation among the ``special'' unipotent classes of $G$. 


\end{abstract}


\section{Introduction} \label{sec0}
Let $\Irr(W)$ be the set of (complex) irreducible characters of a finite
Coxeter group $W$. There is a natural partition $\Irr(W)=\bigsqcup_{\,\cF} 
\Irr(W\mid \cF)$ where $\cF$ runs over the two-sided cells of $W$ in the 
sense of Kazhdan--Lusztig \cite{KaLu}. This partition is an important 
ingredient in the fundamental work of Lusztig \cite{LuBook} on the 
characters of reductive groups over finite fields. Using some standard 
operations in the character ring of $W$ (truncated induction from parabolic 
subgroups, tensoring with the sign character), Lusztig has defined another 
partition of $\Irr(W)$ into so-called ``families''. As shown in 
\cite[Chap.~5]{LuBook} (see also \cite[Chap.~23]{Lusztig03}), these two 
partitions turn out to be the same. The proof relies on deep results from 
algebraic geometry which provide certain ``positivity'' properties of the 
Kazhdan--Lusztig basis \cite{KaLu} of the associated Iwahori--Hecke algebra. 
 
Now, the theory of Kazhdan--Lusztig cells gives rise not only to the 
partition $\Irr(W)=\bigsqcup_{\,\cF}\Irr(W\mid\cF)$, but also to a natural 
partial order $\leq_{\cLR}$ on the pieces in this partition. For example, 
if $W$ is the symmetric group $\fS_n$, then $\Irr(W)$ is parametrized by 
the partitions of $n$, all families are singleton sets, and $\leq_{\cLR}$ 
corresponds to the dominance order on partitions; see \cite{mymurphy} and 
the references there. This is the prototype of a picture which applies to
any finite $W$. 

The main purpose of this paper is to obtain a better understanding of 
the partial order $\leq_{\cLR}$. This will be relevant in a number of 
applications; we just mention, for example, that $\leq_{\cLR}$ is a 
crucial ingredient in defining a ``cellular structure'' (in the sense
of Graham--Lehrer \cite{GrLe}) of the associated Iwahori--Hecke algebra 
\cite{mycell}. Our first main result will show that $\leq_{\cLR}$ can be 
characterised in a purely elementary way in terms of standard operations 
in the character ring of $W$ (induction, truncated induction, tensoring 
with sign), similar in spirit to Lusztig's definition of families. In 
particular, we obtain an efficient algorithm for computing the partial 
order, which can be implemented in {\sf CHEVIE} \cite{chevie}. We 
conjecture that this remains valid in the more general framework of 
Lusztig \cite{Lusztig83}, \cite{Lusztig03} where ``weights'' may be 
attached to the generators of $W$. (We provide both theoretical and 
experimental evidence for this conjecture.)

The main inspiration for this work is a paper by Spaltenstein \cite{Spalt}. 
By pushing the ideas in \cite{Spalt} a little bit further, and combining
them with the above characterisation of $\leq_{\cLR}$, we obtain our second 
main result:

{\em If $W$ is the Weyl group of an algebraic group $G$, then the partial 
order $\leq_{\cLR}$ on the families of $\Irr(W)$ can be interpreted, via 
the Springer correspondence, in terms of the closure relation among 
the ``special'' unipotent classes of $G$}. 

This paper is organised as follows. We recall the basic definitions on
cells and families in Section~\ref{sec1}. Here, we work in the general 
framework of Iwahori--Hecke algebras with unequal parameters, taking into
account ``weight functions'' as in \cite{Lusztig83}, \cite{Lusztig03}. In 
Definition~\ref{mydef} and Conjecture~\ref{mainc}, we propose our 
alternative description of $\leq_{\cLR}$ (in the form of an equivalence). 
In Section~\ref{sec2}, we prove at least one implication in that 
conjectured equivalence in the general case of unequal parameters; see 
Proposition~\ref{prop31}. This is followed by the discussion of some 
examples in which the reverse implication can be seen to hold by 
elementary methods. In Section~\ref{sec3}, we concentrate on the equal
parameter case and complete the proof of Conjecture~\ref{mainc} in that 
case. This allows us to discuss in Section~\ref{sec4} the relation with 
unipotent classes and the work of Spaltenstein \cite{Spalt}.

It would be interesting to understand how our results in Section~\ref{sec4}
are related to work of Bezrukavnikov \cite[Theorem~4]{bezru}. In a 
completely different direction, by work of Brou\'e, Chlouveraki, Kim, 
Malle, Rouquier (see \cite{Chlou}), there is also a notion of ``families'' 
for the irreducible characters of finite complex reflection groups. It 
would be interesting to see if it is possible to define a partial order
on these families as well. (As Jean Michel has pointed out to me, one
cannot simply adopt the definitions in this paper.)

\section{Kazhdan--Lusztig cells and families} \label{sec1}

Let $W$ be a finite Coxeter group, with generating set $S$ and corresponding 
length function $l\colon W \rightarrow \Z_{\geq 0}$. Let $\Gamma$ be an 
abelian group (written additively) and $L \colon W \rightarrow \Gamma$ be a 
weight function, that is, we have $L(ww')=L(w)+L(w')$ whenever $w,w'\in W$ 
are such that $l(ww')=l(w)+l(w')$. Let $F \subseteq \C$ be a splitting field
for $W$ and $A=F[\Gamma]$ be the $F$-vector space with basis $\{v^g \mid 
g\in \Gamma\}$. There is a well-defined ring structure on $A$ such that 
$v^gv^{g'}= v^{g+g'}$ for all $g,g' \in \Gamma$. Let $\bH=\bH_A(W,S,L)$ be 
the corresponding generic Iwahori--Hecke algebra over $A$ with parameters 
$\{v_s \mid s\in S\}$ where $v_s:=v^{L(s)}$ for $s\in S$. This is an 
associative algebra which is free as an $A$-module, with basis $\{T_w\mid 
w \in W\}$. The multiplication is given by the rule
\[ T_sT_w=\left\{\begin{array}{cl} T_{sw} & \quad \mbox{if $l(sw)>l(w)$},\\
T_{sw}+(v_s-v_s^{-1})T_w & \quad \mbox{if $l(sw)<l(w)$},\end{array}
\right.\]
where $s\in S$ and $w\in W$.  See \cite{gepf}, \cite{Lusztig83}, 
\cite{Lusztig03} for further details.

We assume that there exists a total ordering $\leq$ of $\Gamma$ which is 
compatible with the group structure, that is, whenever $g,g',h \in \Gamma$ 
are such that $g\leq g'$, then $g+h\leq g'+h$. This implies that $A$ is an 
integral domain; we denote by $K$ its field of fractions. Throughout this 
paper, we assume that
\[ L(s)\geq 0 \qquad \mbox{for all $s \in S$}.\]
We define $\Gamma_{\geq 0}=\{g\in \Gamma\mid g\geq 0\}$ and denote by
$\Z[\Gamma_{\geq 0}]\subseteq A$ the set of all integral linear combinations
of terms $v^g$ where $g\geq 0$. The notations $\Z[\Gamma_{>0}]$,
$\Z[\Gamma_{\leq 0}]$, $\Z[\Gamma_{<0}]$ have a similar meaning.

\begin{exmp} \label{equalp} Let $\Gamma=\Z$ and $\leq$ be the natural
order. (This is the setting of Lusztig \cite{Lusztig03}.) Then $A$ is 
nothing but the ring of Laurent polynomials over $F$ in the 
indeterminate~$v$. We have $K=F(v)$. If, furthermore, we have $L(s)=1$ for 
all $s \in S$, then we say that we are in the ``equal parameter case''.
\end{exmp}

Returning to the general case, let $\{C_w \mid w \in W\}$ be the 
Kazhdan--Lusztig basis of $\bH$; see \cite{KaLu}, \cite{Lusztig83}, 
\cite{Lusztig03}. The element $C_w$ is characterised by the property that 
(a) it is fixed by a certain ring involution of $\bH$ and (b) it is 
congruent to $T_w$ modulo $\sum_{y\in W} \Z[\Gamma_{>0}] T_y$. (This is the 
original convention used in \cite{KaLu}, \cite{Lusztig83}.) Let 
$\leq_{\cL}$, $\leq_{\cR}$, $\leq_{\cLR}$ be the Kazhdan--Lusztig pre-order 
relations on $W$; for any $w \in W$, we have
\[ \bH C_w \subseteq \sum_{y \leq_{\cL} w} {\Z}[\Gamma]C_y,\qquad
C_w\bH \subseteq \sum_{y \leq_{\cR} w} {\Z}[\Gamma]C_y,\qquad
\bH C_w \bH\subseteq \sum_{y \leq_{\cLR} w} {\Z}[\Gamma]C_y.\]
Let $\sim_{\cL}$, $\sim_{\cR}$, $\sim_{\cLR}$ be the associated equivalence
relations on $W$. Thus, given $x,y \in W$, we have $x \sim_{\cL} y$ if and
only if $x \leq_{\cL} y$ and $y \leq_{\cL} x$. (Similarly for $\sim_{\cR}$ 
and $\sim_{\cLR}$.) The corresponding equivalence classes are called ``left
cells'', ``right cells'' and ``two-sided cells'', respectively. Note that 
all these notions depend on the weight function $L$ and the total ordering of 
$\Gamma$. 

Let $\fC$ be a left cell and set $[\fC]_A:=\fI_\fC/\hat{\fI}_\fC$ where
\begin{align*}
\fI_\fC &=\mbox{$A$-span}\{C_y \mid y \leq_{\cL} w 
\mbox{ for some $w \in \fC$}\}\\
\hat{\fI}_\fC &=\mbox{$A$-span}\{C_y \mid \,y \leq_{\cL} w \mbox{ for 
some $w \in \fC$, but $y \not\in \fC$}\}.
\end{align*}
Since $\fI_\fC$ and $\hat{\fI}_\fC$ are left ideals in $\bH$, the quotient
$[\fC]_A$ is a left $\bH$-module with a canonical $A$-basis indexed by the
elements of $\fC$. Extending scalars from $A$ to $F$ via the $F$-algebra
homomorphism $\theta_1 \colon A \rightarrow F$ sending all $v^g$ to $1$ 
($g \in \Gamma$), we obtain a left $F[W]$-module $[\fC]_1:=F \otimes_A 
[\fC]_A$. We have a direct sum decomposition of left $F[W]$-modules 
\[ F[W] \cong \bigoplus_{\text{$\fC$ left cell in $W$}} [\fC]_1.\]
Now let us denote by $\Irr(W)$ the set of irreducible representations of
$W$ over $F$ (up to isomorphism); recall that $F$ is assumed to be a 
splitting field for $W$. Let $E \in \Irr(W)$. Since we have the above 
direct sum decomposition, there exists a left cell $\fC$ such that $E$ is 
a constituent of $[\fC]_1$; furthermore, all such left cells are contained 
in the same two-sided cell. This two-sided cell, therefore, only depends 
on $E$ and will be denoted by $\cF_E$.  Thus, we obtain a natural surjective 
map
\[ \Irr(W) \rightarrow \{ \mbox{set of two-sided cells of $W$}\}, 
\quad E\mapsto \cF_E.\]
(See Lusztig \cite[5.15]{LuBook} for the equal parameter case; the same 
argument works in general.) It will be useful to introduce the following
notation. Let $X,Y$ be any subsets of $W$. Then we write $X \leq_{\cLR} Y$ 
if $x \leq_{\cLR} y$ for all $x \in X$ and $y \in Y$. 

\begin{defn}{\bf (Lusztig \protect{\cite{LuBook}})} \label{family1} Let 
$E,E' \in \Irr(W)$. We write $E \leq_{\cLR} E'$ if $\cF_E \leq_{\cLR}
\cF_{E'}$. This defines a pre-order relation on $\Irr(W)$. We write 
$E \sim_{\cLR} E'$ if $E \leq_{\cLR} E'$ and $E' \leq_{\cLR} E$ or, 
equivalently, if $\cF_E=\cF_{E'}$. Thus, we obtain a partition
\[\Irr(W)=\bigsqcup_{\cF \text{ two-sided cell}} \Irr(W\mid \cF),\]
where $\Irr(W\mid \cF)$ consists of all $E \in \Irr(W)$ such that 
$\cF_E=\cF$.
\end{defn}

\begin{rem} \label{dual} Let $w_0\in W$ be the longest element and
$\sgn$ be the sign representation of $W$. If $\fC$ is a left cell in
$W$, then so is $\fC w_0$ and we have 
\[ [\fC w_0]_1 \cong [\fC]_1 \otimes \sgn.\]
(See Lusztig \cite[Lemma~5.14]{LuBook} and \cite[Cor.~2.8]{my02}.) 
Furthermore, multiplication by $w_0$ reverses the relations $\leq_{\cL}$,
$\leq_{\cR}$ and $\leq_{\cLR}$; see \cite[Cor.~11.7]{Lusztig03}. It 
follows that, for all $E,E'\in \Irr(W)$, we have: 
\begin{itemize}
\item[(a)] $\cF_{E\otimes \sgn}=\cF_E\,w_0$.
\item[(b)] $E \leq_{\cLR} E'$ if and only if $E' \otimes\sgn\leq_{\cLR}
E \otimes \sgn$.
\end{itemize}
Thus, tensoring with $\sgn$ induces an order-reversing bijection on the
sets $\Irr(W\mid \cF)$.
\end{rem}

In order to describe Lusztig's alternative characterisation of the sets 
$\Irr(W\mid\cF)$, we need to introduce some further notation. Recall that 
$K$ is the field of fractions of $A=F[\Gamma]$. By extension of scalars, 
we obtain a $K$-algebra $\bH_K=K\otimes_A \bH$ which is known to be split 
semisimple; see \cite[9.3.5]{gepf}. Furthermore, by Tits' Deformation 
Theorem, the irreducible representations of $\bH_K$ (up to isomorphism) are 
in bijection with the irreducible representations of $W$; see 
\cite[8.1.7]{gepf}. Given $E \in \Irr(W)$, we denote by $E_v$
the corresponding irreducible representation of $\bH_K$. This is 
uniquely characterised by the following condition:
\[ \theta_1(\mbox{trace}(T_w,E_v))=\mbox{trace}(w,E) \qquad 
\mbox{for all $w \in W$},\]
where $\theta_1 \colon A \rightarrow F$ is as above. Note also that 
$\mbox{trace}(T_w,E_v) \in A$ for all $w\in W$.

\begin{defn}{\bf (Lusztig)} \label{ainv} Given $E \in \Irr(W)$, we define
\[ \ba_E:=\min \{g \in \Gamma_{\geq 0} \mid v^g \mbox{trace}(T_w,E_v) 
\in F[\Gamma_{\geq 0}] \mbox{ for all $w \in W$}\}\] 
Furthermore, we define numbers $c_{w,E} \in F$ by 
\[ \mbox{trace}(T_w,E_v)=c_{w,E} \,v^{-\ba_E} + \mbox{combination of 
terms $v^g$ where $g>-\ba_E$}.\]
\end{defn}

(In the equal parameter case, these definitions were given by Lusztig 
\cite[(5.1.21)]{LuBook}. The same definitions work in general; see also
\cite{my02}). The following result shows that the numbers $c_{w,E}$ can, 
in fact, be used to detect the two-sided cell $\cF_E$.

\begin{lem}[Lusztig] \label{lemlus} We have $\varnothing \neq \{w \in W 
\mid c_{w,E}\neq 0\} \subseteq \cF_E$ for all $E \in \Irr(W)$.
\end{lem}

(See Lusztig \cite[Lemma~5.2]{LuBook} for the equal parameter case;
the same arguments also work in general. For more details in the general
case, see \cite[Prop.~4.7]{my02}.)

Now let $I \subseteq S$ and consider the parabolic subgroup $W_I\subseteq W$
generated by $I$. Then we have a corresponding parabolic subalgebra $\bH_I 
\subseteq \bH$. By extension of scalars from $A$ to $K$, we also have a 
subalgebra $\bH_{K,I}=K \otimes_A \bH_I \subseteq \bH_K$. The above
definitions (i.e., $\ba_E$, $c_{w,E}$, $\ldots$) apply to the irreducible
representations of $W_I$ as well. Denote by $\Ind_I^S$ the induction of 
representations, either from $W_I$ to $W$ or from $\bH_I$ to $\bH$.

\begin{lem}[Lusztig] \label{inda} Let $M \in \Irr(W_I)$.
\begin{itemize}
\item[(a)] If $E \in \Irr(W)$ is a constituent of $\Ind_I^S(M)$, then 
$\ba_E \geq \ba_M$.
\item[(b)] There exists some $E \in \Irr(W)$ which is a constituent of 
$\Ind_I^S(M)$ and such that $\ba_E=\ba_M$.
\end{itemize}
\end{lem}

(See Lusztig \cite{Lusztig79b} in the equal parameter case; the same
arguments work in general. See \cite[Lemma~3.5]{my02} for details.)

We now recall Lusztig's definition of families. Let $E \in \Irr(W)$ and 
$M \in\Irr(W_I)$. We write $M \rightsquigarrow_L E$ if $E$ is a constituent 
of $\Ind_I^S(M)$ and $\ba_E= \ba_M$. 

\begin{defn}{\bf (Lusztig \protect{\cite[4.2]{LuBook}})} \label{family2} 
The partition of $\Irr(W)$ into ``families'' is defined as follows. When 
$W= \{1\}$, there is only one family; it consists of the unit representation 
of $W$. Assume now that $W \neq \{1\}$ and that families have already been 
defined for all proper parabolic subgroups of $W$. Then $E,E' \in \Irr(W)$ 
are said to be in the same family for $\Irr(W)$ if there exists a sequence 
$E= E_0,E_1, \ldots, E_m=E'$ in $\Irr(W)$ such that, for each $i \in \{0,1,
\ldots,m-1\}$, the following condition is satisfied. There exists a subset 
$I_i \subsetneqq S$ and $M_i',M_i'' \in \Irr(W_{I_i})$, where $M_i'$, 
$M_i''$ belong to the same family of $\Irr(W_{I_i})$, such that 
either
\begin{align*}
 M_i' \rightsquigarrow_L E_{i-1} \qquad &\mbox{and} \qquad 
M_i'' \rightsquigarrow_L E_i\\
\intertext{or}
M_i' \rightsquigarrow_L E_{i-1} \otimes \sgn \qquad &\mbox{and} \qquad 
M_i'' \rightsquigarrow_L E_i\otimes \sgn.
\end{align*}
\end{defn}

Note that it is clear from this definition that tensoring with the
sign representation permutes the families. 

We can now state the following remarkable theorem. One of its applications 
is that it facilitates the explicit determination of the partition of 
$\Irr(W)$ in Definition~\ref{family1}; see Lusztig \cite[Chap.~4]{LuBook}.

\begin{thm}[Barbasch--Vogan, Lusztig \protect{\cite[5.25]{LuBook}}] 
\label{lusthm1} Assume that $W$ is a finite Weyl group and that we are in 
the equal parameter case. Let $E,E' \in \Irr(W)$. Then $E \sim_{\cLR} E'$ 
(see Definition~\ref{family1}) if and only if $E,E'$ belong to the same 
family (see Definition~\ref{family2}).
\end{thm}

\begin{rem} \label{lusthm1a} The ``if'' part of the above result is 
proved by elementary methods; see \cite[Chap.~5]{LuBook}. Our 
Proposition~\ref{prop31} below provides a new proof for this ``if'' part, 
which also works in the general multi-parameter case. The proof of the 
``only if'' part in \cite{LuBook} relies on deep results from the theory 
of primitive ideals in enveloping algebras (which also explains the 
restriction to Weyl groups). An alternative approach is provided by 
\cite[23.3]{Lusztig03} and \cite{myert05} where it is shown that the 
above theorem holds for any finite $W$ and any weight function 
$L \colon W\rightarrow\Gamma$, assuming that Lusztig's conjectures 
{\bf P1}--{\bf P15} in \cite[14.2]{Lusztig03} are satisfied. This is 
known to be true for all finite Coxeter groups in the equal parameter 
case (see the comments on the proof of Theorem~\ref{a123} below); it is 
also true for a number of situations involving unequal parameters. For a 
summary of the present state of knowledge, see \cite[\S 5]{myisom} and 
the references there.
\end{rem}

Our aim is to find an alternative description of the pre-order $\leq_{\cLR}$
on $\Irr(W)$, in the spirit of Lusztig's definition of families. The
following definition is inspired by Spaltenstein \cite{Spalt}.

\begin{defn} \label{mydef} We define a relation $\preceq$ on $\Irr(W)$ 
inductively as follows. If $W=\{1\}$, then $\Irr(W)$ only consists of the 
unit representation and this is related to itself. Now assume that $W 
\neq\{1\}$ and that $\preceq$ has already been defined for all proper
parabolic subgroups of $W$. Let $E,E'\in\Irr(W)$. Then we write $E\preceq 
E'$ if there is a sequence $E=E_0,E_1,\ldots, E_m=E'$ in $\Irr(W)$ such 
that, for each $i \in \{0,1,\ldots,m-1\}$, the following condition is 
satisfied. There exists a subset $I_i \subsetneqq S$ and $M_i',M_i''\in 
\Irr(W_{I_i})$, where $M_i' \preceq M_i''$ within $\Irr(W_{I_i})$, such 
that either 
\begin{align*}
E_{i-1} \mbox{ is a constituent of } \Ind_{I_i}^S(M_i') \qquad &\mbox{and} 
\qquad M_i'' \rightsquigarrow_L E_i\\\intertext{or}
E_i \otimes \sgn \mbox{ is a constituent of } \Ind_{I_i}^S(M_i')\qquad &
\mbox{and} \qquad M_i'' \rightsquigarrow_L E_{i-1}\otimes \sgn.
\end{align*}
\end{defn}

We note that, as in \cite[4.2]{LuBook}, it is enough to require that, 
in the above definition, we have $|I_i|=|S|-1$ for all $i$ (that is,
each $W_{I_i}$ is a maximal parabolic subgroup).

\begin{rem} \label{mydef1} Let $E,E'\in \Irr(W)$. It is clear from the
above definition that we have the following implications:
\begin{itemize}
\item[(a)] If $E,E'$ belong to the same family then $E \preceq 
E'$ and $E'\preceq E$.
\item[(b)] If $E\preceq E'$, then we also have $ E'\otimes \sgn \preceq 
E \otimes \sgn$.
\end{itemize}
The reverse implication in (a) does not seem to follow easily from
the definitions. In Proposition~\ref{prop31a}, we will establish that 
reverse implication in the equal parameter case; the general multi-parameter
case requires further work and will be dealt with in \cite[Cor.~6.2]{klord2}.
\end{rem}

By analogy with Theorem~\ref{lusthm1}, we would now like to state the 
following:

\begin{conj} \label{mainc} Let $E,E' \in \Irr(W)$. Then $E \leq_{\cLR} E'$
(see Definition~\ref{family1}) if and only if $E \preceq E'$ (see
Definition~\ref{mydef}). 
\end{conj}

In Section~\ref{sec2}, we will prove the ``if'' part of the conjecture
by a general argument (for any weight function $L \colon W \rightarrow
\Gamma$ as above). In particular, as already announced in
Remark~\ref{lusthm1a}, this will provide a new, completely elementary 
proof of the ``if'' part of Theorem~\ref{lusthm1}. We also verify in some 
examples that the reverse implications hold. In Section~\ref{sec3}, we 
will prove the ``only if'' part of the conjecture by a general argument,
assuming that we are in the equal parameter case. 

\section{Two-sided cells and induced representations} \label{sec2}

We keep the setting of the previous section, where $W$ is a finite Coxeter
group and $L \colon W \rightarrow \Gamma$ is any weight function such that 
$L(s)\geq 0$ for all $s \in S$.

Given a subset $I \subseteq S$, let $W_I$ be the corresponding parabolic
subgroup of $W$ and $X_I$ be the set of distinguished left coset 
representatives of $W_I$ in $W$. Thus, we have a bijection $X_I \times 
W_I \rightarrow W$, $(d,w)\mapsto dw$, where $l(dw)=l(d)+l(w)$; 
see \cite[\S 2.1]{gepf}. In the following discussion, we shall make 
frequent use of the main result of \cite{myind}, concerning the
induction of cells from $W_I$ to $W$.

\begin{lem} \label{lem10} Let $x\in W$ and write $x=dw$ where $d \in 
X_I$ and $w\in W_I$. Let $E \in \Irr(W)$ be a constituent of $[\fC]_1$ 
where $\fC$ is the left cell of $W$ which contains~$x$; in particular, 
$x \in \cF_E$. Then there exists some $M \in \Irr(W_I)$ such that $w 
\in \cF_M$ and $E$ is a constituent of $\Ind_I^S(M)$. 
\end{lem}

\begin{proof} Let $\fC'$ be the left cell in $W_I$ which contains~$w$.  
Then, by \cite[Theorem~1]{myind}, we have $\fC \subseteq X_I \fC'$; 
furthermore, by \cite[Lemma~5.2]{myert05}, $[\fC]_1$ is a direct summand
of $\Ind_I^S([\fC']_1)$. Hence, since $E \in \Irr(W)$ is a constituent 
of $[\fC]_1$, there exists some $M \in \Irr(W_I)$ such that $M$ is a
constituent of $[\fC']_1$ and $E$ is a constituent of $\Ind_I^S(M)$. 
Since $w \in \fC'$, we also have $w \in \cF_M$, as required. 
\end{proof}

Recall that, for any subsets $X,Y$ of $W$, we write $X \leq_{\cLR} Y$ if 
$x \leq_{\cLR} y$ for all $x \in X$ and $y \in Y$. 

\begin{lem} \label{lem11} Let $E \in\Irr(W)$ and $M \in \Irr(W_I)$ be 
such that $E$ is a constituent of $\Ind_I^S(M)$. Then we have 
$\cF_E \leq_{\cLR} \cF_M$.
\end{lem}

\begin{proof} Let $\fC'$ be a left cell in $W_I$ such that $M$ is 
a constituent of $[\fC']_1$. As in the above proof, by 
\cite[Theorem~1]{myind}, we have a partition $X_I \fC'=\bigsqcup_{i=1}^m 
\fC_i$ where $\fC_1,\ldots, \fC_m$ are left cells of $W$. Furthermore, 
by \cite[Lemma~5.2]{myert05}, we have $\Ind_I^S([\fC']_1)=
\bigoplus_{i=1}^m [\fC_i]_1$. Hence, since $E$ is a constituent of 
$\Ind_I^S(M)$, there exists some $i$ such that $E$ is a constituent of 
$[\fC_i]_1$. Let $\fC:=\fC_i$. Now note that $l(xw)=l(x)+ l(w)$ for all 
$x \in X_I$ and $w\in W_I$. This length condition implies that $xw 
\leq_{\cL} w$ for all $x \in X_I$ and $w \in W_I$; see 
\cite[Theorem~6.6]{Lusztig03}. Hence, we have $w\leq_{\cL} w'$ for all 
$w \in \fC$ and $w'\in \fC'$. Since $\fC'\subseteq \cF_M$ and $\fC 
\subseteq \cF_E$, this implies that $\cF_E \leq_{\cLR} \cF_M$, as
required.
\end{proof}

A special case of the following result appeared in 
\cite[Lemma~3.6]{mymurphy}.

\begin{lem} \label{lem12} Let $E \in\Irr(W)$ and $M \in \Irr(W_I)$ be 
such that $M \rightsquigarrow_L E$. Then we have $\cF_M \subseteq \cF_E$.
\end{lem}

\begin{proof} The algebra $\bH$ is symmetric, with trace form $\tau\colon
\bH\rightarrow A$ given by $\tau(T_1)=1$ and $\tau(T_w)=0$ for $1 \neq w 
\in  W$. The sets $\{T_w \mid w \in W\}$ and $\{T_{w^{-1}}\mid w \in W\}$ 
form a pair of dual bases. Hence we have the following orthogonality 
relations:
\[ \sum_{w \in W} \mbox{trace}(T_w,E_v) \,\mbox{trace}(T_{w^{-1}},
E_v')=\left\{\begin{array}{cl} (\dim E)\,\bc_E & \quad 
\mbox{if $E \cong E'$}, \\ 0 & \quad \mbox{otherwise};
\end{array}\right.\]
see \cite[8.1.8]{gepf}. Here, $0 \neq \bc_E \in A$ and, as observed
by Lusztig, we have 
\[ \bc_E=f_E\, v^{-2\ba_E}+ \mbox{combination of terms $v^g$ where 
$g>-2\ba_E$},\]
where $f_E$ is a strictly positive real number; see \cite[3.3]{my02}.
The same definitions apply, of course, to the parabolic subalgebra
$\bH_I$. Now consider the element
\[ e_M:=\sum_{w \in W_I} \mbox{trace}(T_w,M_v)\, T_{w^{-1}}
\in \bH_{K,I}.\]
We shall evaluate $\mbox{trace}(e_M,E_v)$ in two ways.  On the 
one hand, given $E'\in \Irr(W)$, let us denote by $d(E',M)$ the
multiplicity of $E'$ as a constituent of $\Ind_I^S(M)$. By Frobenius 
reciprocity and the compatibility with specialisations in 
\cite[9.1.9]{gepf}, this implies that 
\[ \mbox{trace}(h,E_v)=\sum_{M' \in \Irr(W_I)} d(E,M')\, 
\mbox{trace}(h, M_v')\qquad\mbox{for all $h \in \bH_{K,I}$}.\]
Using the orthogonality relations for the irreducible representations of
$\bH_{K,I}$, we conclude that 
\begin{align*} 
\mbox{trace}(e_M,E_v) &=  \sum_{M' \in \Irr(W_I)} d(E,M')\,
\mbox{trace}(e_M,M_v')\\ &= \sum_{M' \in \Irr(W_I)} d(E,M')\, 
\sum_{w \in W_I} \mbox{trace}(T_w,M_v)\, \mbox{trace}(T_{w^{-1}},M_v')\\
&= d(E,M)\,(\dim M)\, \bc_M.
\end{align*}
Consequently, we have 
\[ v^{2\ba_M}\,\mbox{trace}(e_M,E_v)=d(E,M)\, (\dim M)\,
f_M+ \mbox{``higher terms''},\]
where ``higher terms'' means an $F$-linear combination of terms 
$v^g$ where $g \in \Gamma_{>0}$. On the other hand, recalling 
Definition~\ref{ainv} and taking into account our assumption 
$\ba_M=\ba_E$, we obtain
\begin{align*}
v^{2\ba_M}\,\mbox{trace}(e_M,E_v) &=  \sum_{w \in W_I} 
\bigl(v^{\ba_M}\,\mbox{trace} (T_w,M_v)\bigr)\,\bigl(v^{\ba_E}\,
\mbox{trace}(T_{w^{-1}}, E_v)\bigr)\\
&=\Bigl(\sum_{w \in W_I} c_{w,M}\, c_{w^{-1},E}\Bigr)+
\mbox{``higher terms''}.
\end{align*}
Comparing the two expressions, we deduce that 
\begin{equation*}
\sum_{w \in W_I} c_{w,M}\, c_{w^{-1},E}=d(E,M)\, (\dim M)\, f_M. \tag{$*$}
\end{equation*}
Now the right hand side of ($*$) is non-zero since $d(E,M) \neq 0$ by
assumption. Hence, there exists some $w \in W_I$ such that $c_{w,M} 
\neq 0$ and $c_{w^{-1},E} \neq 0$. By \cite[Cor.~8.2.6]{gepf}, we have 
$\mbox{trace}(T_w, E_v)=\mbox{trace}(T_{w^{-1}},E_v)$.
So we also have $c_{w,E}=c_{w^{-1},E}\neq 0$. By Lemma~\ref{lemlus}, this 
implies $w \in \cF_M \cap \cF_E$ and, hence, $\cF_M \subseteq \cF_E$.
\end{proof}

\begin{prop} \label{prop31} Let $E,E'\in \Irr(W)$. If $E \preceq E'$, then
$E \leq_{\cLR} E'$. In particular, if $E,E'$ belong to the same family, 
then $E \sim_{\cLR} E'$.
\end{prop}

\begin{proof} If $W=\{1\}$, there is nothing to prove. Now assume that
$W\neq \{1\}$ and that the assertion has already been proved for all
proper parabolic subgroups of $W$. It is now sufficient to consider an 
elementary step in Definition~\ref{mydef}. That is, we can assume that 
there is a subset $I \subsetneqq S$ and $M',M'' \in \Irr(W_I)$, where 
$M' \preceq M''$ within $\Irr(W_I)$, such that one of the following two 
conditions holds.
\begin{itemize}
\item[(I)] $E$ is a constituent of $\Ind_I^S(M')$ and $M'' 
\rightsquigarrow_L E'$.
\item[(II)] $E'\otimes \sgn$ is a constituent of $\Ind_I^S(M')$ and 
$M'' \rightsquigarrow_L E\otimes \sgn$.
\end{itemize}
If (I) holds, then $\cF_E \leq_{\cLR} \cF_{M'}$ and $\cF_{M''} \subseteq 
\cF_{E'}$ by Lemmas~\ref{lem11} and~\ref{lem12}. Since $M'\preceq M''$, 
we already know that $M'\leq_{\cLR} M''$ and, hence, $\cF_{M'}\leq_{\cLR} 
\cF_{M''}$ (with respect to $W_I$). But then we also have 
$\cF_{M'} \leq_{\cLR} \cF_{M''}$ with respect to $W$ and, hence, $\cF_E 
\leq_{\cLR} \cF_{E'}$, as required. 

On the other hand, if (II) holds, then a completely similar argument shows 
that $\cF_{E'\otimes \sgn} \leq_{\cLR} \cF_{E\otimes \sgn}$. But, by
Remark~\ref{dual}, we have $\cF_{E \otimes \sgn}=\cF_Ew_0$ and 
$\cF_{E' \otimes \sgn}=\cF_{E'}w_0$. Furthermore, multiplication
with $w_0$ reverses the relation $\leq_{\cLR}$. Hence, we have 
$\cF_{E} \leq_{\cLR} \cF_{E'}$, as required. 

Finally, if $E,E'$ belong to the same family, then Remark~\ref{mydef1}
immediately shows that $E \preceq E'$, $E' \preceq E$ and, hence, 
$E \sim_{\cLR} E'$. 
\end{proof}

\begin{exmp} \label{exph4} Let $(W,S)$ be of type $H_4$. Here, we are
automatically in the equal parameter case. There are $34$ irreducible 
representations in $\Irr(W)$ and they are partitioned into $13$ families;
see Alvis--Lusztig \cite{AlLu82}. Using {\sf CHEVIE} \cite{chevie}, one 
easily determines the relation $\preceq$. It turns out that we obtain a 
``linear'' order such that, for all $E,E'\in \Irr(W)$, we have: 
\begin{itemize}
\item[(a)] $E \preceq E'$ is and only if $\ba_{E'} \leq \ba_{E}$.
\item[(b)] $E,E'$ belong to the same family if and only if $\ba_{E}=
\ba_{E'}$.
\end{itemize}
On the other hand, Alvis \cite{Alvis87} has determined the two-sided 
cells of $W$; there are precisely $13$ of them. Hence, by 
Proposition~\ref{prop31}, we have $E \sim_{\cLR} E'$ if and only if 
$E,E'$ belong to the same family. Furthermore, since $\preceq$ already 
induces a linear order on families, it follows that $E \preceq E'$ if 
and only if $E \leq_{\cLR} E'$.  Thus, Conjecture~\ref{mainc} holds 
in this case.

Similar remarks apply to $(W,S)$ of type $H_3$ and $I_2(m)$ (with equal 
or unequal parameters in the latter case): In all these cases, one easily
checks that $\preceq$ is a linear order satisfying (a), (b) and, hence, 
Conjecture~\ref{mainc} holds. (See the summary of the relevant results
on cells and families in \cite[\S 7]{myert05}.)
\end{exmp}

\begin{exmp} \label{expf4} Let $(W,S)$ be of type $F_4$, with generators
and diagram given by 
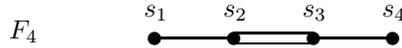
\begin{center}
\begin{picture}(200,20)
\put( 10, 5){$F_4$}
\put( 61,13){$s_1$}
\put( 91,13){$s_2$}
\put(121,13){$s_3$}
\put(151,13){$s_4$}
\put( 65, 5){\circle*{5}}
\put( 95, 5){\circle*{5}}
\put(125, 5){\circle*{5}}
\put(155, 5){\circle*{5}}
\put( 65, 5){\line(1,0){30}}
\put( 95, 7){\line(1,0){30}}
\put( 95, 3){\line(1,0){30}}
\put(125, 5){\line(1,0){30}}
\end{picture}
\end{center}
Let $\Gamma=\Z$ and $L$ be a weight function which is specified by two 
positive integers $a:=L(s_1)=L(s_2)>0$ and $b:=L(s_3)=L(s_4)>0$. Taking 
into account the symmetry of the diagram, one may assume that $a\leq b$. 
There are $25$ irreducible representations in $\Irr(W)$. The relation 
$\leq_{\cLR}$ on $\Irr(W)$ has been determined in all cases in 
\cite{compf4}. It turns out that there are only four essentially different 
cases: $b=a$, $b=2a$, $2a>b>a$ or $b>2a$; see Table~1 in
\cite[p.~362]{compf4}.

It is verified in \cite{compf4} that $E \sim_{\cLR} E'$ if and only if 
$E,E'$ belong to the same family. Using {\sf CHEVIE} \cite{chevie}, one
easily determines the relation $\preceq$. By inspection, one finds
that Conjecture~\ref{mainc} holds in all cases. One also finds that:
\begin{itemize}
\item[(a)] If $E \preceq E'$ then $\ba_{E'} \leq \ba_{E}$.
\item[(b)] If $E \preceq E'$ and $\ba_E=\ba_{E'}$, then $E,E'$ belong to 
the same family.
\end{itemize}
This example provides strong evidence for the validity of 
Conjecture~\ref{mainc} in the general case of unequal parameters.
\end{exmp}

\begin{exmp} \label{expbn} Let $(W,S)$ be of type $B_n$, with  
generators and diagram given by 
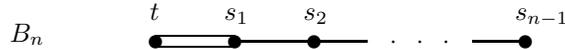
\begin{center}
\begin{picture}(250,20)
\put( 10, 5){$B_n$}
\put( 63,13){$t$}
\put( 91,13){$s_1$}
\put(121,13){$s_2$}
\put(201,13){$s_{n{-}1}$}
\put( 65, 5){\circle*{5}}
\put( 95, 5){\circle*{5}}
\put(125, 5){\circle*{5}}
\put( 65, 7){\line(1,0){30}}
\put( 65, 3){\line(1,0){30}}
\put( 95, 5){\line(1,0){30}}
\put(125, 5){\line(1,0){20}}
\put(155,5){\circle*{1}}
\put(165,5){\circle*{1}}
\put(175,5){\circle*{1}}
\put(185, 5){\line(1,0){20}}
\put(205, 5){\circle*{5}}
\end{picture}
\end{center}
We have $\Irr(W)=\{E^\lambda \mid \lambda \in \Lambda\}$ where $\Lambda$ 
is the set of all pairs of partitions of total size $n$. For example, 
the unit, sign and reflection representation are labelled by $((n),
\varnothing)$, $(\varnothing,(1^n))$ and $((n-1), (1))$, respectively;
see \cite[\S 5.5]{gepf}. Let $\Gamma=\Z$. Then a weight function $L$ is 
specified by two integers $b:=L(t)\geq 0$ and $a=L(s_i)\geq 0$ for 
$1\leq i\leq n-1$. For a conjectural description of the partial order 
$\leq_{\cLR}$ on two-sided cells, see \cite[Remark~1.2]{bgil}. 

Here is a specific example in the case of unequal parameters, where we 
assume that $b>(n-1)a>0$. This is the ``asymptotic'' case originally 
studied by Bonnaf\'e and Iancu \cite{BI}, \cite{BI2}. By 
Proposition~\ref{prop31} and \cite[Prop.~5.4]{geia06}, we have 
\[ E^\lambda \preceq E^\mu \qquad \Rightarrow \qquad E^\lambda 
\leq_{\cLR} E^\mu \qquad \Rightarrow \qquad \lambda\trianglelefteq \mu\]
where $\trianglelefteq$ denotes the dominance order on pairs of partitions.
In order to prove the reverse implications, it will be enough to show 
that $\lambda \trianglelefteq \mu \Rightarrow E^\lambda \preceq E^\mu$. 
Thus, we are reduced to a purely combinatorial problem. This, and a full 
description of $\preceq$ for all choices of the parameters $a,b$, will 
be discussed in \cite{klord2}.
\end{exmp}

\section{The equal parameter case} \label{sec3}

Throughout this section, we assume that $\Gamma=\Z$ and $L(s)=1$ for 
all $s \in S$. Our aim is to show that, in this setting, 
Conjecture~\ref{mainc} holds. For this purpose, we have to rely on 
some deep properties of the relations $\leq_{\cL}$, $\leq_{\cR}$, 
$\leq_{\cLR}$ which are stated in Theorem~\ref{a123} below. These in
turn are established by using certain ``positivity'' properties of the 
Kazhdan--Lusztig basis of $\bH$ which are only available in the equal 
parameter case; see Lusztig \cite[Chap.~16]{Lusztig03} and the references 
there (as far as finite Weyl groups are concerned) and DuCloux \cite{fokko} 
(as far as types $H_3$, $H_4$, $I_2(m)$ are concerned). 

\begin{thm} \label{a123} In the equal parameter case, the following hold.
\begin{itemize}
\item[(a)] {\rm (Lusztig \cite{Lusztig03})} If $E,E'\in \Irr(W)$ are such 
that $E\leq_{\cLR} E'$, then $\ba_{E'}\leq \ba_E$. In particular, if $E 
\sim_{\cLR} E'$, then $\ba_E=\ba_{E'}$. 
\item[(b)] {\rm (Lusztig \cite{Lusztig03})} If $E,E'\in \Irr(W)$ are such 
that $E\leq_{\cLR} E'$ and $\ba_{E'}=\ba_E$, then $E \sim_{\cLR} E'$.
\item[(c)] {\rm (Lusztig--Xi \cite{LusXi})} Let $x,y\in W$ be such that $x
\leq_{\cLR} y$. Then there exists some $z \in W$ such that $x \leq_{\cL} 
z$ and $z \sim_{\cR} y$.
\end{itemize}
\end{thm}

{\em Comments on the proof}. Using the ``positivity'' properties mentioned
above, Lusztig shows in \cite[Chap.~16]{Lusztig03} that the conjectural
properties {\bf P1}--{\bf P15} in \cite[14.2]{Lusztig03} hold for $\bH$.
Then (a) and (b) are a combination of {\bf P4}, {\bf P11} and 
\cite[Prop.~20.6]{Lusztig03}. The statement in (c) is due to Lusztig--Xi 
\cite[\S 3]{LusXi}. Note that, in \cite{LusXi}, this result is stated for 
affine Weyl groups; but the same proof works when $W$ is finite. Indeed, 
besides general properties of the relations $\leq_{\cL}$, $\leq_{\cR}$, 
$\leq_{\cLR}$, the ingredients needed in the proof are listed in 
\cite[2.2, 2.3, 2.5]{LusXi}. Now, the references for these properties 
cover also the case of finite Coxeter groups; the above-mentioned 
``positivity'' properties are required here, too. An additional reference 
for \cite[2.2(h)]{LusXi} (which is attributed to Springer, unpublished) is 
provided by \cite[1.3]{nxi}. \hfill $\Box$

\begin{rem} \label{rem123} By Lusztig's conjectures in 
\cite[14.2]{Lusztig03}, one can expect that (a) and (b) remain valid
in the general case of unequal parameters. The proof of (c) seems to 
require more than just using the conjectural properties {\bf P1}--{\bf P15}
in \cite[14.2]{Lusztig03}. It is not clear (at least not to me) if one can 
expect (c) to hold in the general case of unequal parameters.
\end{rem}

As a first application of Theorem~\ref{a123}(a), we obtain the following 
converse to Lemma~\ref{lem12}.

\begin{lem} \label{lem10a} Let $I \subseteq S$. Let $E \in \Irr(W)$ 
and $M \in \Irr(W_I)$ be such that $\cF_M\subseteq \cF_E$ and $E$ 
is a constituent of $\Ind_I^S(M)$. Then $M \rightsquigarrow_L E$.
\end{lem}

\begin{proof} By Lemma~\ref{inda}(b), there exists some $E'\in \Irr(W)$ 
which is a constituent of $\Ind_I^S(M)$ and such that $\ba_{E'}=\ba_M$. 
By Lemma~\ref{lem12}, we have $\cF_M \subseteq \cF_{E'}$. Thus, we have 
$\cF_M \subseteq \cF_E \cap \cF_{E'}$ and so $\cF_E=\cF_{E'}$. Using 
Theorem~\ref{a123}(a), we conclude that $\ba_E=\ba_{E'}=\ba_M$, as required. 
\end{proof}

Next recall from Remark~\ref{mydef1} that, if $E,E'\in \Irr(W)$ belong
to the same family, then $E \preceq E'$ and $E' \preceq E$. Now we can also
prove the reverse implication.

\begin{prop} \label{prop31a} Let $E,E'\in \Irr(W)$ be such that
$E \preceq E'$. Then $\ba_{E'} \leq \ba_E$. Furthermore, if $E\preceq E'$
and $E'\preceq E$, then $\ba_E=\ba_{E'}$ and $E,E'$ belong to the same
family of $\Irr(W)$.
\end{prop}

\begin{proof} By Proposition~\ref{prop31}, we have $E \leq_{\cLR} E'$.
So Theorem~\ref{a123}(a) implies that $\ba_{E'}\leq \ba_E$. Now assume
that $E \preceq E'$ and $E'\preceq E$. Then, clearly, $\ba_E=\ba_{E'}$.

We now show by an inductive argument that, if $E \preceq E'$ and
$\ba_E=\ba_{E'}$, then $E, E'$ belong to the same family. If $W=\{1\}$, 
there is nothing to prove. Now assume that $W \neq \{1\}$ and that the 
assertion has already been proved for all proper parabolic subgroups of 
$W$. As in the proof of Proposition~\ref{prop31}, it is sufficient to 
consider an elementary step in Definition~\ref{mydef}. That is, we can 
assume that there is a subset $I \subsetneqq S$ and $M',M''\in\Irr(W_I)$, 
where $M' \preceq M''$ within $\Irr(W_I)$, such that one of the following 
two conditions holds.
\begin{itemize}
\item[(I)] $E$ is a constituent of $\Ind_I^S(M')$ and $M'' 
\rightsquigarrow_L E'$.
\item[(II)] $E' \otimes \sgn$ is a constituent of $\Ind_I^S(M')$ and $M'' 
\rightsquigarrow_L E \otimes \sgn$.
\end{itemize}
First of all, since $M '\preceq M''$, we already know that $\ba_{M''}
\leq \ba_{M'}$.

Now, if (I) holds, then $\ba_E\geq \ba_{M'} \geq \ba_{M''}=\ba_{E'}$. 
Since $\ba_E=\ba_{E'}$, we conclude that $\ba_{M'}=\ba_{M''}$. Hence, 
by induction, $M',M''$ belong to the same family of $\Irr(W_I)$. 
Furthermore, since $\ba_{E}=\ba_{M'}$, we have $M' \rightsquigarrow_L E$. 
Thus, the first set of conditions in Definition~\ref{family2} is satisfied
and so $E,E'$ belong to the same family of $\Irr(W)$.

On the other hand, if (II) holds, then $\ba_{E' \otimes \sgn}\geq 
\ba_{M'}\geq \ba_{M''}=\ba_{E \otimes \sgn}$. Assume, if possible, that
$\ba_{E' \otimes \sgn}>\ba_{E \otimes \sgn}$. Then $E \otimes \sgn
\not\sim_{\cLR} E' \otimes \sgn$ by Theorem~\ref{a123}(a). Consequently,
we also have $E \not\sim_{\cLR} E'$ by Remark~\ref{dual}. Since $E 
\leq_{\cLR} E'$ and $\ba_E=\ba_{E'}$, this contradicts 
Theorem~\ref{a123}(b). Hence, we must have $\ba_{E' \otimes\sgn}=
\ba_{E\otimes \sgn}$. Now we can argue as above and conclude that
the second set of conditions in Definition~\ref{family2} is satisfied. 
Hence, $E,E'$ belong to the same family of $\Irr(W)$.
\end{proof}

(Note that the above proof only requires (a) and (b) in Theorem~\ref{a123}.)

\begin{rem} \label{remgeia} In \cite[Cor.~6.1]{klord2} we will show that
Proposition~\ref{prop31a} remains valid in the general multi-parameter 
case. The proof relies on a case-by-case argument and a detailed study of 
the relation $\preceq$ in type $B_n$.
\end{rem}

Besides the above-mentioned ``positivity'' properties, another 
distinguished feature of the equal parameter case is the existence of 
``special'' irreducible representations. (As discussed in 
\cite[Example~4.11]{compf4}, one cannot expect the existence of 
representations with similar properties in the general case of unequal 
parameters.) Given $E \in \Irr(W)$, let $\bb_E$ be the smallest 
$i \geq 0$ such that $E$ is a constituent of the $i$-th symmetric power 
of the natural reflection representation of $W$. It is an empirical 
observation that we always have $\ba_E \leq \bb_E$; following Lusztig 
\cite{Lusztig79b}, we say that $E$ is ``special'' if $\ba_E=\bb_E$. Let 
\[\cS(W):=\{E \in \Irr(W) \mid E \mbox{ special} \}.\]

\begin{thm}[Lusztig \protect{\cite[4.14]{LuBook}}] \label{Lspecial} Each 
family of $\Irr(W)$ (see Definition~\ref{family2}) contains a unique 
$E \in \cS(W)$. 
\end{thm}

(See also \cite[\S 6.5]{gepf} where non-crystallographic Coxeter groups
are included from the outset in the discussion.)

\begin{thm}[Lusztig \protect{\cite[5.25]{LuBook}}] \label{Lspecial1} 
Let $\fC$ be a left cell and let $E \in \cS(W)$. If $\fC \subseteq 
\cF_E$, then $E$ occurs with multiplicity~$1$ in $[\fC]_1$.
\end{thm}

(Alternative proofs are provided by \cite{Lusztig86}, \cite{myert05}; 
these references also cover the cases where $W$ is of type $H_3$, $H_4$ 
or $I_2(m)$.) 

\begin{rem} \label{Lspecial2a} Let $I \subsetneqq S$ and let $\cS(W_I)$
denote the set of all $M \in \Irr(W_I)$ which are special (with respect to
$W_I$). Let $M \in \cS(W_I)$. Then it is known (see \cite{Lusztig79b}) that 
there is a unique $E \in \cS(W)$ such that $\ba_E=\ba_M$ and $\Ind_I^S(M)$ 
equals $E$ plus a sum of irreducible representations $E' \in \Irr(W)$ such 
that $\ba_{E'}> \ba_E$; in particular, we have $M \rightsquigarrow_L E$. 
Let us write $E=j_{I}^S(M)$ in this case. 
\end{rem}

We define $\cS^\circ(W)$ to be the set of all $j_I^S(M)$ where 
$I \subsetneqq S$ and $M \in \cS(W_I)$.  With this definition, we can 
now state the following result of Spaltenstein which will be a further
key ingredient in our argument.

\begin{lem}[Cf.\ Spaltenstein \protect{\cite{Spalt}}] \label{lemspalt}
Let $E \in \cS(W)$ be such that $E \not\in \cS^\circ(W)$. Then 
$\ba_{E \otimes\sgn} <\ba_E$.
\end{lem}

\begin{proof} By standard reduction arguments, it is enough to prove
this in the case where $(W,S)$ is irreducible. If $W$ is of type $H_3$, 
$H_4$ or $I_2(m)$, the assertion is easily checked by an explicit 
computation and {\sf CHEVIE} \cite{chevie}. One could also check the
assertion for finite Weyl groups in this way, using the explicit knowledge
of $\cS(W)$ and of the invariants $\ba_E$ from \cite{Lusztig79b}. However,
a related verification has already been done by Spaltenstein 
\cite[\S 5]{Spalt}. Thus, all we need to do is to see how the setting in 
\cite[\S 5]{Spalt} translates to our setting here.

So now assume that $W$ is a finite Weyl group. Let $G$ be a simple algebraic 
group (over $\C$ or over $\overline{\F}_p$ where $p$ is a large prime) with 
Weyl group $W$. Using the Springer correspondence (see \cite{Spr}, 
\cite{LuIC}), we can naturally associate with every $E \in \Irr(W)$ a pair 
consisting of a unipotent class of $G$, which we denote by $O_E$, and a 
$G$-equivariant irreducible local system on $O_E$. By \cite[13.1.1]{LuBook}, 
we have
\[ \dim \cB_u=\ba_E \qquad \mbox{for $E \in \cS(W)$},\]
where $\cB_u$ denotes the variety of Borel subgroups containing an element
$u \in O_E$. 

Now Spaltenstein \cite[\S 5]{Spalt} has shown that, if $E \in \cS(W)$ 
and $E \not\in \cS^\circ(W)$, then $O_E$ is strictly contained in the 
Zariski closure of $O_{\bar{E}}$ where $\bar{E}$ is the unique special 
representation of $W$ in the same family as $E \otimes\sgn$. In particular, 
we have $\dim \cB_{\bar{u}}<\dim \cB_u$ where $u \in O_E$ and $\bar{u}
\in O_{\bar{E}}$. Hence, we also have $\ba_{\bar{E}}< \ba_E$. Finally, 
by Proposition~\ref{prop31} and Theorem~\ref{a123}(a), we have 
$\ba_{\bar{E}}=\ba_{E \otimes\sgn}$.
\end{proof}

Given a two-sided cell $\cF$ in $W$, we denote by $\ba(\cF)$ the common 
value of $\ba_E$ where $E \in \Irr(W)$ is such that $\cF_E=\cF$; see 
Theorem~\ref{a123}(a). With this convention, we can now state the
following version of Lemma~\ref{lemspalt} which does not refer to 
``special'' representations in $\Irr(W)$. (One may conjecture that this 
remains true in the general case of unequal parameters.)

\begin{cor} \label{lemspalt1} Let $\cF$ be a two-sided cell in $W$ such 
that $\cF \cap W_I=\varnothing$ for all proper subsets $I\subsetneqq S$.
Then $\ba(\cF w_0)<\ba(\cF)$.
\end{cor}

\begin{proof} By Proposition~\ref{prop31} and Theorem~\ref{Lspecial}, 
there exists some $E \in \cS(W)$ such that $\cF_E=\cF$. Assume, if
possible, that there exists some $I \subsetneqq S$ and $M \in \cS(W_I)$
such that $E=j_I^S(M)$.  In particular, this would mean that $E$ is a
constituent of $\Ind_I^S(M)$ and $\ba_M=\ba_E$. Hence, by 
Lemma~\ref{lem12}, we would have $\cF_M \subseteq \cF_E=\cF$ and so 
$\cF \cap W_I \neq \varnothing$, a contradiction. Thus, we have
$E \not\in \cS^\circ(W)$. Now Lemma~\ref{lemspalt} implies that
$\ba_{E \otimes\sgn}<\ba_E$. 

By Remark~\ref{dual}, we have $\cF_{E \otimes\sgn}=\cF_E w_0$. Hence,
we have $\ba_E=\ba(\cF_E)$ and $\ba_{E \otimes\sgn}=\ba(\cF_E w_0)$. This 
yields $\ba(\cF w_0)<\ba(\cF)$, as required.
\end{proof}

\begin{thm} \label{mainthm} Recall our standing assumption that we are in 
the equal parameter case. Then Conjecture~\ref{mainc} holds.  
\end{thm}

\begin{proof} The ``if'' part is already proved in Proposition~\ref{prop31}.
To prove the ``only if'' part, we use an inductive argument. If $W=\{1\}$, 
there is nothing to prove. Now assume that $W\neq \{1\}$ and that the 
``only if'' part has already been proved for all proper parabolic subgroups 
$W$. Let $E,E' \in \Irr(W)$ be such that $E \leq_{\cLR} E'$. We must show 
that $E \preceq E'$. Since $E \leq_{\cLR} E'$, we have $\cF_E \leq_{\cLR} 
\cF_{E'}$. We claim that one of the following two conditions is satisfied:
\begin{itemize}
\item[(I)] $\cF_{E'} \cap W_I\neq \varnothing$ for some $I \subsetneqq S$.
\item[(II)] $\cF_Ew_0\cap W_I\neq \varnothing$ for some $I \subsetneqq S$.
\end{itemize}
To prove this, we use an argument due to Spaltenstein \cite{Spalt}. Assume, 
if possible, that $\cF_{E'} \cap W_I=\varnothing$ and $\cF_Ew_0 \cap W_I=
\varnothing$ for all $I\subsetneqq S$. By Corollary~\ref{lemspalt1}, this 
implies that $\ba(\cF_{E'}w_0)<\ba(\cF_{E'})$ and $\ba(\cF_E)<\ba(\cF_E
w_0)$. Furthermore, since $\cF_E \leq_{\cLR} \cF_{E'}$, we have 
$\ba(\cF_{E'}) \leq \ba(\cF_E)$ by Theorem~\ref{a123}(a). Thus, we 
conclude that $\ba(\cF_{E'}w_0)<\ba(\cF_{E}w_0)$. On the other hand, since 
$\cF_E \leq_{\cLR} \cF_{E'}$, we also have $\cF_{E'}w_0 \leq_{\cLR} 
\cF_Ew_0$ (see Remark~\ref{dual}). So, Theorem~\ref{a123}(a) implies that
$\ba(\cF_Ew_0)\leq \ba(\cF_{E'}w_0)$, and we have reached a contradiction. 
Thus, (I) or (II) holds, as claimed. 

Now let us first assume that (I) holds. Let $E_0$ be the unique special 
representation in the same family as $E$ and $E_0'$ be the unique special 
representation in the same family as $E'$; see Theorem~\ref{Lspecial}. Then 
$E \preceq E_0$ and $E_0' \preceq E'$ by Remark~\ref{mydef1}(a). Hence, it 
will be enough to show that $E_0 \preceq E_0'$. Note that, by 
Proposition~\ref{prop31}, we have $\cF_E=\cF_{E_0}$ and $\cF_{E'}=
\cF_{E_0'}$.

Let $y \in \cF_{E'} \cap W_I$. Then we claim that there exists some 
$x\in\cF_E$ such that $x\leq_{\cL} y$. This is seen as follows. Recall 
from Remark~\ref{dual} that multiplication by the longest element $w_0
\in W$ reverses the relations $\leq_{\cL}$, $\leq_{\cR}$ and 
$\leq_{\cLR}$. Now take any element $x' \in \cF_E$. Since $\cF_E 
\leq_{\cLR} \cF_{E'}$, we have $x'\leq_{\cLR} y$. Then $yw_0\leq_{\cLR} 
x'w_0$ and so, by Theorem~\ref{a123}(c), there exists some $z \in W$ 
such that $yw_0\leq_{\cL} z$ and $z \sim_{\cR} x'w_0$. In particular, 
$z \in \cF_E w_0$ and so $x:=zw_0 \in \cF_E$. Since $yw_0 \leq_{\cL} z=
xw_0$, we now deduce that $x\leq_{\cL} y$, as required. 

Let us write $x=dw$ where $d \in X_I$ and $w \in W_I$, as in 
Lemma~\ref{lem10}. Thus, $x=dw \leq_{\cL} y$ where $y \in W_I$. Then, by
relation ($\dagger$) in \cite[\S 4]{myind}, we have $w \leq_{\cLR,I} 
y$ where the subscript $I$ indicates that this relation is with respect 
to $W_I$. 

Let $\fC$ be the left cell in $W$ which contains~$x$. Then $E_0$ is 
a constituent of $[\fC]_1$; see Theorem~\ref{Lspecial1}. By 
Lemma~\ref{lem10}, there exists some $M\in \Irr(W_I)$ such that 
$w \in \cF_M$ and $E_0$ is a constituent of $\Ind_I^S(M)$. Similarly,
let $\fC'$ be the left cell in $W$ which contains~$y$; now $E_0'$ is 
a constituent of $[\fC']_1$. Again, there exists some $M' \in \Irr(W_I)$ 
such that $y\in \cF_{M'}$ and $E_0'$ is a constituent of $\Ind_I^S(M')$. 
Furthermore, since $y \in \cF_{M'} \cap \cF_{E_0'}$, we must have 
$\cF_{M'} \subseteq \cF_{E_0'}$. So we can now conclude that $M' 
\rightsquigarrow_L E_0'$; see Lemma~\ref{lem10a}. Since $w\leq_{\cLR,I} 
y$, we have $\cF_{M} \leq_{\cLR,I} \cF_{M'}$ and so $M \leq_{\cLR,I} M'$. 
By our inductive hypothesis, we deduce that $M \preceq M'$ within 
$\Irr(W_I)$. Thus, the first set of conditions in Definition~\ref{mydef} 
is satisfied. Hence, we have $E_0 \preceq E_0'$ and so $E\preceq E'$. 
This completes the proof in the case where (I) holds. 

Finally, assume that (II) holds. Then we can argue as follows. By 
Remark~\ref{dual}, we have $\cF_{E \otimes\sgn}=\cF_{E}w_0$ and
$\cF_{E' \otimes\sgn}=\cF_{E'}w_0$. In particular, (II) is equivalent to 
$\cF_{E\otimes\sgn} \cap W_I \neq \varnothing$. Furthermore, since 
$\cF_E \leq_{\cLR} \cF_{E'}$, we have $\cF_{E'\otimes\sgn}\leq_{\cLR} 
\cF_{E \otimes\sgn}$. We can now apply the same argument as above and 
conclude that $E'\otimes\sgn \preceq E \otimes\sgn$. Then 
Remark~\ref{mydef1}(b) shows that we also have $E \preceq E'$, as required.
\end{proof}

\section{Unipotent classes and two-sided cells} \label{sec4}

We continue to assume that we are in the equal parameter case. In
addition, we now assume that $W$ is the Weyl group of a connected 
reductive algebraic group $G$ (over $\C$ or over $\overline{\F}_p$ 
where $p$ is a large prime). By the Springer correspondence (see 
\cite{Spr}, \cite{LuIC}), we can naturally associate with every $E \in 
\Irr(W)$ a pair consisting of a unipotent class of $G$, which we denote 
by $O_E$, and a $G$-equivariant irreducible local system on $O_E$. Thus, 
we obtain a map 
\[ \Irr(W) \rightarrow \{\mbox{set of unipotent classes in $G$}\},
\quad E \mapsto O_E.\]
(The local system on $O_E$ will not play a role for our purposes here.) 

\begin{defn}{\bf (Lusztig)} \label{defsp} A unipotent class $O$ of 
$G$ is called ``special'' if $O=O_E$ where $E\in \cS(W)$. The map 
$E \mapsto O_E$ gives a bijection between $\cS(W)$ and the set of special 
unipotent classes in $G$.
\end{defn}

\begin{rem} \label{defsp1} Let $\cF$ be a two-sided cell in $W$
and consider the collection of unipotent classes
\[\cC(\cF):=\{O_E\mid E\in\Irr(W) \mbox{ such that }\cF_E= \cF\}.\]
By Theorems~\ref{lusthm1} and~\ref{Lspecial}, there exists a unique 
$E_0 \in \cS(W)$ such that $\cF_{E_0}=\cF$; in particular, $O_{E_0}\in
\cC(\cF)$. Then it is known that 
\[O\subseteq\overline{O}_{E_0}\qquad\mbox{for all $O\in\cC(\cF)$};\]
see \cite[Prop.~2.2]{GeMa2}. (Here, and below, $\overline{X}$ denotes the 
Zariski closure in $G$ for any subset $X \subseteq G$.) Thus, the 
special unipotent class $O_{E_0}$ can be characterized as the unique 
unipotent class in $\cC(\cF)$ which is maximal with respect to the 
Zariski closure relation.
\end{rem}

Let $\cU_G$ be the unipotent variety of $G$. Let $O$ be a special 
unipotent class. The corresponding ``special piece'' in $\cU_G$ is 
defined to be the set of all elements in $\overline{O}$ which are not 
contained in $\overline{O}'$ where $O'$ is any special unipotent class 
such that $\overline{O}'\subsetneqq\overline{O}$. By Spaltenstein 
\cite{Spa} and Lusztig \cite{Lspec}, the special pieces form a partition 
of $\cU_G$. Note that every special piece is a union of a special 
unipotent class (which is open dense in the special piece) and of a
certain number (possibly zero) of non-special unipotent classes. 

We can now associate with every two-sided cell in $W$ a special piece 
in $\cU_G$, as follows. Let $\cF$ be a two-sided cell in $W$. As already
noted above, there exists a unique $E_0 \in \cS(W)$ such that $\cF_{E_0}=
\cF$. Let $O_{E_0}$ be the corresponding special unipotent class and 
$\cO_{\cF}$ be the unique special piece in $\cU_G$ containing $O_{E_0}$. 
Thus, we obtain a canonical bijection (see also Lusztig
\cite[Theorem~0.2]{Lspec}):
\[\{\mbox{set of two-sided cells of $W$}\} 
\stackrel{\text{1{-}1}}{\longrightarrow} \{\mbox{set of special 
pieces in $\cU_G$}\}, \quad \cF \mapsto \cO_{\cF}.\]
As remarked in \cite[\S 14]{L12}, this map is part of Lusztig's bijection 
\cite{aff4} between the set of two-sided cells in an associated affine 
Weyl group and the set of {\em all} unipotent classes of $G$.

Corollary~\ref{fthm} below gives an interpretation of the order relation
$\leq_{\cLR}$ on the two-sided cells of $W$ in terms of the closure 
relation among the special pieces in $\cU_G$. This will heavily rely on 
Theorem~\ref{mainthm} and on the following result.

\begin{thm}[Spaltenstein \protect{\cite{Spa}, \cite{Spalt}}] 
\label{spaltthm} Let $E,E'\in \cS(W)$. Then we have 
\[ E \preceq_s E' \quad \Leftrightarrow \quad O_{E} \subseteq 
\overline{O}_{E'} \quad \Leftrightarrow \quad  O_{\bar{E}'}
\subseteq \overline{O}_{\bar{E}}.\]
\end{thm}

Here, we have used the following notation. Given $E \in \cS(W)$, we denote
by $\bar{E} \in \cS(W)$ the unique special representation in the same 
family as $E \otimes\sgn$. (Thus, we obtain an involution $E \mapsto 
\bar{E}$ on $\cS(W)$.) Furthermore, the relation $\preceq_s$ on $\cS(W)$ 
is defined inductively as follows. If $W=\{1\}$, then $\cS(W)$ only 
consists of the unit representation and this is related to itself. Now 
assume that $W \neq\{1\}$ and that $\preceq_s$ has already been defined for 
all proper parabolic subgroups of $W$. Let $E,E'\in\cS(W)$. Then we 
write $E\preceq_s E'$ if there exists a subset $I \subsetneqq S$ and 
$M',M''\in \cS(W_I)$, where $M' \preceq_s M''$ within $\cS(W_I)$, 
such that either 
\begin{align*}
E \mbox{ is a constituent of } \Ind_I^S(M') \qquad &\mbox{and} \qquad 
M'' \rightsquigarrow_L E'\\\intertext{or}
\bar{E}' \mbox{ is a constituent of } \Ind_I^S(M')\qquad &\mbox{and} 
\qquad M'' \rightsquigarrow_L \bar{E}.
\end{align*}
Note the formal similarity in the definitions of $\preceq_s$ and the 
relation $\preceq$ considered in Section~\ref{sec1}.  More precisely,
we have:

\begin{lem} \label{lem4} Let $E,E' \in \cS(W)$ be such that $E\preceq_s 
E'$. Then we also have $E \preceq E'$.
\end{lem}

\begin{proof} We proceed by an inductive argument. If $W=\{1\}$, there 
is a nothing to prove. Now assume that $W=\{1\}$ and that the assertion 
has already been proved for all proper parabolic subgroups of $W$. By 
the definition of $\preceq_s$, there exists a subset $I \subsetneqq S$ 
and $M',M''\in \cS(W_I)$, where $M' \preceq_s M''$ within $\cS(W_I)$, such 
that one of the folllowing conditions is satisfied:
\begin{itemize}
\item[(I)] $E$ is a constituent of $\Ind_I^S(M')$ and $M'' 
\rightsquigarrow_L E'$.
\item[(II)] $\bar{E}'$ is a constituent of $\Ind_I^S(M')$ and $M'' 
\rightsquigarrow_L \bar{E}$.
\end{itemize}
By our inductive hypothesis, we already know that $M' \preceq M''$ within 
$\Irr(W_I)$. Consequently, if (I) holds, then the first set of conditions 
in Definition~\ref{mydef} is satisfied and so $E \preceq E'$. Now assume
that (II) holds. Then we obtain that $\bar{E}' \preceq \bar{E}$. By the 
definition of $\bar{E}$, $\bar{E}'$ and Remark~\ref{mydef1}(a), we have 
$\bar{E} \preceq E\otimes\sgn$, $E' \otimes\sgn \preceq \bar{E'}$ and
so $E'\otimes\sgn\preceq E\otimes\sgn$. Hence, Remark~\ref{mydef1}(b) 
implies that $E \preceq E'$, as required. 
\end{proof}

\begin{lem} \label{resform} Let $P\subseteq G$ be a parabolic subgroup 
of $G$, with unipotent radical $U_P$ and Levi complement $L$ such that 
$L$ has Weyl group $W_I\subseteq W$ where $I \subseteq S$. Let $E \in
\Irr(W)$ and $O_E$ be the corresponding unipotent class in $G$; let
$M \in \Irr(W_I)$ and $O_M$ be the corresponding unipotent class in $L$.
\begin{itemize}
\item[(a)] Assume that $E$ is a constituent of $\Ind_I^S(M)$. Then 
$O_E\cap U_PO_M \neq \varnothing$. 
\item[(b)] Assume that $E$ is special and $M \rightsquigarrow_L E$. Then 
$M$ is special and $U_P O_M \subseteq \overline{O}_E$.
\end{itemize}
\end{lem}

\begin{proof} (a) Springer's restriction formula \cite[Theorem~4.4]{Spr} 
(see also Lusztig \cite[Theorem~8.3]{LuIC}) expresses the multiplicity of 
$E$ as a constituent of $\Ind_I^S(M)$ in geometric terms, using the variety 
\[\fX_{u,u'}(P):=\{x \in G \mid x^{-1}ux \in u'U_P\}, \quad
\mbox{where $u \in O_E$ and $u' \in O_M$}.\]
In particular, the assumption that $E$ is a constituent of $\Ind_I^S(M)$ 
implies that $\fX_{u,u'}(P)$ must be non-empty. Thus, we have $O_E\cap 
U_P O_M\neq\varnothing$, as required.

(b) We check that $O_E$ is induced from $O_M$ in the sense of
Lusztig--Spaltenstein \cite{LuSp}. To begin with, since $E$ is special, 
the unipotent class $O_E$ has property (B) in \cite[\S 3]{LuSp}; see the 
remark at the end of \cite[\S 2]{Lusztig79b}, or 
\cite[Theorem~6.5.13(c)]{gepf}. On the other hand, since $M 
\rightsquigarrow_L E$, the representation $M$ must also be special.
(This follows, for example, from \cite[\S 5.2 and \S 6.5]{gepf}.) In 
particular, property (B) holds for $O_M$, too. Then \cite[Theorem~3.5]{LuSp}
shows that $O_E$ is induced from $O_M$, that is, $O_E$ is the unique 
unipotent class in $G$ such that $O_E\cap U_P O_M$ is dense in $U_PO_M$. 
Hence, $U_P O_M$ must be contained in the closure of $O_E$, as desired.
\end{proof}

We can now state the promised geometric interpretation of $\leq_{\cLR}$. 

\begin{cor} \label{fthm} Let $\cF$, $\cF'$ be two-sided cells in $W$. Then 
we have $\cF \leq_{\cLR} \cF'$ if and only if $\cO_{\cF} \subseteq 
\overline{\cO}_{\cF'}$.
\end{cor}

\begin{proof} First assume that $\cF \leq_{\cLR} \cF'$. The following 
argument for proving $\cO_{\cF}\subseteq \overline{\cO}_{\cF'}$ is
inspired by the discussion in \cite[\S 2]{Spalt}. If $W=\{1\}$, there is 
nothing to prove. Now assume that $W\neq \{1\}$ and that the assertion has 
already been proved for all proper parabolic subgroups of $W$. As in the 
proof of Theorem~\ref{mainthm}, one of the following two conditions must be 
satisfied:
\begin{itemize}
\item[(I)] $\cF'\cap W_I\neq\varnothing$ for some $I \subsetneqq S$.
\item[(II)] $\cF w_0\cap W_I \neq\varnothing$ for some $I \subsetneqq S$.
\end{itemize}
Assume first that (I) holds. Let $E,E'\in \cS(W)$ be such that $\cF=\cF_E$ 
and $\cF'=\cF_{E'}$. Then we must show that $O_E\subseteq\overline{O}_{E'}$.
As in the proof of Theorem~\ref{mainthm}, since $E,E'$ are special and 
$E\leq_{\cLR} E'$, there exist $M,M'\in\Irr(W_I)$, where $M\leq_{\cLR,I}
M'$ with respect to $W_I$, such that $E$ is a constituent of $\Ind_I^S(M)$ 
and $M' \rightsquigarrow_L E'$. Now let $P\subseteq G$ be a parabolic 
subgroup of $G$, with unipotent radical $U_P$ and Levi complement $L$ 
such that $L$ has Weyl group $W_I$. Applying Lemma~\ref{resform}, we 
conclude that $M'$ is special and that we have the following relations 
among the associated unipotent classes: 
\begin{equation*}
O_{E} \cap U_PO_M \neq \varnothing \qquad \mbox{and}\qquad U_PO_{M'} 
\subseteq \overline{O}_{E'}.\tag{$*$}
\end{equation*} 
Let $M_0 \in \cS(W_I)$ be the unique special representation in the same 
family as $M$ (with respect to $W_I$). Then $M_0 \sim_{\cLR,I} M$ (see
Proposition~\ref{prop31}) and so $M_0\leq_{\cLR,I} M'$. Hence, applying 
our inductive hypothesis, we can conclude that $O_{M_0} \subseteq 
\overline{O}_{M'}$ (within $L$). Furthermore, since $M, M_0$ belong to 
the same family, we have $O_M \subseteq \overline{O}_{M_0}$; see 
Remark~\ref{defsp1}. Thus, we have reached the conclusion that $O_{M} 
\subseteq \overline{O}_{M'}$. This certainly implies that $U_PO_M$ is 
contained in the closure of $U_PO_{M'}$. Combining this with ($*$), it 
follows that $O_{E} \subseteq \overline{O}_{E'}$, as required. 

Now assume that (II) holds. Then the same argument shows that 
$O_{\bar{E}'} \subseteq \overline{O}_{\bar{E}}$; note that, by
Proposition~\ref{prop31} and Remark~\ref{dual}, we have $\cF_{\bar{E}_0}
=\cF_{E_0\otimes\sgn}=\cF_{E_0}w_0$ for every $E_0 \in \cS(W)$. But, by 
the second equivalence in Theorem~\ref{spaltthm}, we then also have that 
$O_{E} \subseteq \overline{O}_{E'}$, as required.

Conversely, assume that $\cO_{\cF} \subseteq \overline{\cO}_{\cF'}$.
Let again $E,E'\in \cS(W)$ be such that $\cF=\cF_E$ and $\cF'=\cF_{E'}$. 
Then the assumption certainly implies that $O_{E} \subseteq 
\overline{O}_{E'}$. So the first equivalence in Theorem~\ref{spaltthm}
shows that $E \preceq_s E'$. By Lemma~\ref{lem4} and 
Proposition~\ref{prop31}, this implies $E \preceq E'$ and $E \leq_{\cLR}
E'$, as required.
\end{proof}

\begin{rem} \label{frem} The closure relation among the special unipotent 
classes in $G$, and the order-reversing bijection $O_{E}\mapsto 
O_{\bar{E}}$ ($E \in \cS(W)$), are explicitly known; see Carter 
\cite[\S 13.2]{Carter2}, Spaltenstein \cite{Spa}. Hence, by the above 
result, we also have an explicit description of the partial order 
$\leq_{\cLR}$ on the families of $\Irr(W)$.

On the other hand, the advantage of Theorem~\ref{mainthm} is that it
provides a purely elementary description of $\leq_{\cLR}$ in terms of
the relation $\preceq$, independently of the theory of algebraic groups. 
Moreover, the equivalence between $\leq_{\cLR}$ and $\preceq$ applies 
to more general situations where no geometric interpretation is 
available; see the examples in Section~\ref{sec2}.
\end{rem}

{\it Note added in proof.} After the submission of this paper, I 
learned that the statement of Corollary~\ref{fthm} already appeared 
as Proposition~2.23 in an article by Barbasch and Vogan, Annals of Math. 
{\bf 121} (1985), 41--110. The details of the proof of the ``if'' 
part are omitted there, and the proof of the ``only if'' part is 
different from the one given here.

 

\begin{thebibliography}{131} 
 
\bibitem{Alvis87}
{\sc D.~Alvis}, The left cells of the {C}oxeter group of type ${H}_4$. J.
Algebra {\bf 107} (1987), 160--168.

\bibitem{AlLu82}
{\sc D.~Alvis and G.~Lusztig}, The representations and generic degrees of the
 {H}ecke algebra of type ${H}_4$. J. Reine Angew. Math. \textbf{336} (1982),
 201--212; correction, {\em ibid.} \textbf{449} (1994), 217--218.

\bibitem{bgil}
{\sc C.~Bonnaf\'e, M.~Geck, L.~Iancu and T.~Lam},  On domino insertion
and Kazhdan--Lusztig cells in type $B_n$. {\it In:} Representation theory of
algebraic groups and quantum groups (Nagoya, 2006; eds. A.~Gyoja et al.),
Progress in Math. {\bf 284}, Birkh\"auser, 2010.

\bibitem{bezru}
{\sc R. Bezrukavnikov}, Perverse sheaves on affine flags and nilpotent 
cone of the Langlands dual group. Israel J. Math. {\bf 170} (2009), 185--206. 

\bibitem{BI2}
{\sc C.~Bonnaf\'e},  Two-sided cells in type $B$ in the asymptotic
case. J. Algebra {\bf 304} (2006), 216--236.

\bibitem{BI}
{\sc C.~Bonnaf\'e and L.~Iancu}, Left cells in type $B_n$ with unequal
parameters. Represent. Theory {\bf 7} (2003), 587--609.

\bibitem{Carter2}
{\sc R.~W. Carter}, Finite groups of Lie type: Conjugacy classes and 
complex characters. Wiley, New York, 1985; reprinted 1993 as Wiley 
Classics Library Edition.

\bibitem{Chlou}
{\sc M. Chlouveraki}, Blocks and families for cyclotomic Hecke algebras. 
Lecture Notes in Mathematics 1981, Springer, Berlin, 2009.

\bibitem{fokko}
{\sc F.~DuCloux}, Positivity results for the Hecke algebras of
noncrystallographic finite Coxeter group. J. Algebra {\bf 303} (2006),
731--741.

\bibitem{my02}
{\sc M.~Geck}, Constructible characters, leading coefficients and left cells
for finite Coxeter groups with unequal parameters. Represent. Theory {\bf 6}
(2002), 1--30 (electronic).

\bibitem{myind}
{\sc M.~Geck}, On the induction of Kazhdan--Lusztig cells. Bull. London
Math. Soc. {\bf 35} (2003), 608--614.

\bibitem{compf4}
{\sc M.~Geck}, Computing Kazhdan--Lusztig cells for unequal parameters.
J. Algebra {\bf 281} (2004), 342--365; section "Computational Algebra".

\bibitem{myert05}
{\sc M.~Geck}, Left cells and constructible representations. Represent.
Theory {\bf 9} (2005), 385--416; Erratum, {\it ibid.} {\bf 11} (2007), 
172--173.

\bibitem{mymurphy}
{\sc M.~Geck}, Kazhdan--Lusztig cells and the Murphy basis. Proc. London
Math. Soc. {\bf 93} (2006), 635--665.

\bibitem{mycell}
{\sc M. Geck}, Hecke algebras of finite type are cellular. Invent. Math.
{\bf 169} (2007), 501--517.

\bibitem{myisom}
{\sc M.~Geck}, On Iwahori--Hecke algebras with unequal parameters and
Lusztig's isomorphism theorem. Jacques Tits special issue, Pure Appl.
Math. Q. (2010), {\it to appear}; preprint at  {\tt arXiv:0711.2522}.

\bibitem{chevie}
{\sc M. Geck, G.~Hi{\ss}, F.~L\"ubeck, G.~Malle and G.~Pfeiffer},
{\sf CHEVIE}-A system for computing and processing generic character tables
for finite groups of Lie type, Weyl groups and Hecke algebras. Appl.
Algebra Engrg. Comm. Comput. {\bf 7} (1996), 175--210.

\bibitem{geia06}
{\sc M. Geck and L.~Iancu}, Lusztig's $a$-function in type $B_n$ in
the asymptotic case. Special issue celebrating the $60$th birthday of
George Lusztig, Nagoya J. Math. {\bf 182} (2006), 199--240.

\bibitem{klord2}
{\sc M. Geck and L.~Iancu}, On the Kazhdan--Lusztig order on families 
in type $B_n$, {\em in preparation}.

\bibitem{GeMa2} 
{\sc M.~Geck and G.~Malle}, On the existence of a unipotent support for 
the irreducible characters  of finite groups of Lie type. Trans.\ Amer.\
Math.\  Soc. {\bf 352} (2000), 429--456.

\bibitem{gepf}
{\sc M. Geck and G. Pfeiffer}, Characters of finite Coxeter groups and
Iwahori--Hecke algebras. London Math. Soc. Monographs, New Series {\bf 21},
Oxford University Press, New York 2000. xvi+446 pp.

\bibitem{GrLe}
{\sc J.~J.~Graham and G.~I.~Lehrer}, Cellular algebras.  Invent.\ Math. 
{\bf 123} (1996), 1--34.

\bibitem{KaLu}
{\sc D.~A.~Kazhdan and G.~Lusztig}, Representations of {C}oxeter groups and
{H}ecke algebras. Invent. Math. {\bf 53} (1979), 165--184.

\bibitem{Lusztig79b}
{\sc G.~Lusztig}, A class of irreducible representations of a finite {W}eyl
group. Indag. Math. \textbf{41} (1979), 323--335.

\bibitem{Lusztig83}
{\sc G.~Lusztig}, Left cells in {W}eyl groups. {\em In:} Lie Group
Representations, I (eds R.~L. R.~Herb and J.~Rosenberg), Lecture Notes
in Mathematics 1024,  Springer, Berlin, 1983, pp.~99--111.

\bibitem{LuBook}
{\sc G.~Lusztig}, Characters of reductive groups over a finite field.
Annals Math.\ Studies, vol. 107, Princeton University Press, 1984.

\bibitem{LuIC}
{\sc G.~Lusztig}, Intersection cohomology complexes on a reductive group.
Invent. Math. {\bf 75}, 205--272 (1984).

\bibitem{Lusztig86}
{\sc G.~Lusztig}, Sur les cellules gauches des groupes de Weyl. C. R. 
Acad.  Sci. Paris \textbf{302} (1986), 5--8.

\bibitem{aff4}
{\sc G.~Lusztig}, Cells in affine Weyl groups IV. J. Fac. Sci. Tokyo U. (IA)
{\bf 36} (1989), 297--328. 

\bibitem{Lspec}
{\sc G.~Lusztig}, Notes on unipotent classes. Asian J. Math. {\bf 1} 
(1997), 194--207.

\bibitem{Lusztig03}
{\sc G.~Lusztig}, Hecke algebras with unequal parameters. CRM Monographs
Ser.~{\bf 18}, Amer. Math. Soc., Providence, RI, 2003.

\bibitem{L12}
{\sc G.~Lusztig}, Twelve bridges from a reductive group to its Langlands 
dual. Contemp.\ Math.\ {\bf 478} (2009), 125--143. 

\bibitem{LuSp}
{\sc G.~Lusztig and N.~Spaltenstein}, Induced unipotent classes. 
J. London Math. Soc. {\bf 19} (1979), 41--52.

\bibitem{LusXi}
{\sc G.~Lusztig and N.~Xi}, Canonical left cells in affine Weyl groups.
Advances in Math. {\bf 72} (1988), 284--288.

\bibitem{Spa}
{\sc N.~Spaltenstein}, Classes unipotentes et sous-groupes de Borel.
Lecture Notes in Mathematics 946, Springer, 1982.

\bibitem{Spalt}
{\sc N.~Spaltenstein}, A property of special representations of
Weyl groups. J. Reine Angew. Math. {\bf 343} (1983), 212--220.

\bibitem{Spr}
{\sc T.~A.~Springer}, A construction of representations of Weyl groups.
Invent. Math. {\bf 44} (1978), 279--293.

\bibitem{nxi}
{\sc N. Xi}, The leading coefficient of certain Kazhdan--Lusztig 
polynomials of the permutation group $S_n$.  J. Algebra {\bf 285} 
(2005), 136--145.
\end{thebibliography}
\end{document}